\documentclass[11pt]{amsart}
\usepackage{geometry}                % See geometry.pdf to learn the layout options. There are lots.
\usepackage{graphicx}
\usepackage{amssymb}
\usepackage{epstopdf}
\usepackage{lpic}
\usepackage{amssymb,amsmath}
\usepackage{pifont}
\usepackage{relsize}
\usepackage{color}
\usepackage[usenames,dvipsnames,svgnames,table]{xcolor}
\usepackage{float}
\usepackage{rotating}
\usepackage{verbatim}
\usepackage{letltxmacro}
\usepackage[normalem]{ulem}
\usepackage{bbm}
\usepackage{mathtools}
\usepackage{amscd}
\usepackage{mathrsfs}
\allowdisplaybreaks[4]
\usepackage{lineno}
\usepackage{tikz}
\usepackage{hyperref}
\usepackage{graphicx} % use this line instead of the above to suppress graphics in draft copies
\usepackage{indentfirst}

\usepackage{amsthm} % facilities for theorem-like environments
\usepackage{amsfonts}
\usepackage{geometry}

\usepackage{glossaries}

\makeglossaries

\newglossaryentry{naturals}
{
	name = {\ensuremath{\mathbb{N}}},
	description = {Nonnegative integers. Thus an $\mathbb{N}$-graded ring has a degree zero piece},
	sort = N
}

\newglossaryentry{[n]}
{
	name = {\ensuremath{[n]}},
	description = {The set $\{1,\dots,n\}$},
	sort = n
}

\newglossaryentry{spec}
{
	name = {\ensuremath{\operatorname{Spec}}},
	description = {The prime spectrum of a ring, i.e. the set of prime ideals, topologized by the Zariski topology},
	sort = Spec
}

\newglossaryentry{coprod}
{
	name = {\ensuremath{\coprod}},
	description = {The coproduct in the category of sets, i.e., the disjoint union},
	sort = coprod,
	symbol = {\ensuremath{\coprod}}
}

\newglossaryentry{height}
{
	name = {\ensuremath{\operatorname{ht}}},
	description = {The height of a prime ideal, i.e. the maximum length of a chain of prime ideals descending from the given one},
	sort = ht
}

\newglossaryentry{partitions}
{
	name = {\ensuremath{\mathscr{P}}},
	description = {The monoid of partitions with at most $n$ parts, describing the shapes of monomials in $R = A[x_1,\dots,x_n]$},
	sort = P
}

\newglossaryentry{conjparts}
{
	name = {\ensuremath{\overline{\mathscr{P}}}},
	description = {The monoid of partitions with parts of size at most $n$, describing the fine grades of elements in the Stanley-Reisner ring $S = A[B_n\setminus\emptyset]$},
	sort = P'
}

\newglossaryentry{Bn}
{
	name = {\ensuremath{B_n}},
	description = {The boolean algebra of subsets of $[n]$, seen as a poset},
	sort = Bn
}

\newglossaryentry{garsia}
{
	name = {\ensuremath{\mathscr{G}}},
	description = {The Garsia map, defined in \ref{def:garsiamap}},
	sort = garsia
}

\newglossaryentry{Grr}
{
	name = {\ensuremath{G_{rr}}},
	description = {The subgroup of a transformation group $G$ that is generated by its reflections and rotations. In the main case, $G$ is a permutation group, and this is the subgroup generated by transpositions, double transpositions, and 3-cycles},
	sort = Grr
}

\renewcommand{\emptyset}{\ensuremath{\varnothing}}

\newcommand{\C}{\ensuremath{\mathbb{C}}}

\newcommand{\Z}{\ensuremath{\mathbb{Z}}}

\newcommand{\F}{\ensuremath{\mathbb{F}}}

\newcommand{\Proj}{\ensuremath{\mathbb{P}}}

\newcommand{\Aut}{\operatorname{Aut}}

\newcommand{\Ind}{\operatorname{Ind}}
\newcommand{\Res}{\operatorname{Res}}

\theoremstyle{plain}
\newtheorem{thm}{Theorem}[section]
\newtheorem{prop}[thm]{Proposition}
\newtheorem{lemma}[thm]{Lemma}
\newtheorem{cor}[thm]{Corollary}

\theoremstyle{definition}
\newtheorem{definition}[thm]{Definition}
\newtheorem{remark}[thm]{Remark}
\newtheorem{observation}[thm]{Observation}
\newtheorem{question}[thm]{Question}

\newtheorem{notation}[thm]{Convention}

\title{Purely noncommuting groups}
\author{Ben Blum-Smith and Fedor Bogomolov}
\begin{document}
\maketitle

\begin{abstract}
In this paper we define and investigate a class of groups characterized by a representation-theoretic property we call {\em purely noncommuting} or {\em PNC}. This property guarantees that the group has an action on a smooth projective variety with mild quotient singularities. It has intrinsic group-theoretic interest as well. The main results are as follows. (i) All supersolvable groups are PNC. (ii) No nonabelian finite simple groups are PNC. (iii) A metabelian group is guaranteed to be PNC if its commutator subgroup's cyclic prime-power-order factors are all distinct, but not in general. We also give a criterion guaranteeing a group is PNC if its nonabelian subgroups are all large, in a suitable sense, and investigate the PNC property for permutations.
\end{abstract}

%In the first section, we define and investigate this class. In the second, we describe the results of a preliminary search for non-PNC groups that still have the desired action.

\section{Introduction}

A fundamental fact in linear algebra is that any pair of diagonalizable commuting matrices shares a full basis of eigenvectors. If a pair of matrices fails to commute, then they may still share some eigenvectors, although not a full basis. In this case, they act by restriction on the subspace spanned by the common eigenvectors, and their actions on this subspace do commute. Thus one may see the sharing of eigenvectors as a kind of {\em partial commuting}. If two diagonalizable matrices do not share any eigenvectors, they {\em noncommute purely}.

In the representation theory of a finite group on an algebraically closed field of characteristic zero, group elements always act as diagonalizable transformations. In this context, if two elements of an abstract group do commute, then in every representation they will be forced to share a full basis of eigenvectors. But if they do not commute abstractly, it may still be the case that in every concrete representation of this group on a vector space, they are forced to share some eigenvectors, i.e. commute partially. This prompts us to ask: given a finite group $G$, and two elements $x,y\in G$ that do not commute, is it possible to find a representation of $G$ in which this fact is expressed in an unadulterated way, i.e. their abstract failure to commute is realized in a pair of transformations that do not share an eigenvector?

This question motivates the following definitions. Throughout, $G$ is a finite group. All vector spaces are over $\C$. %If we have a representation $\rho:G\rightarrow GL(V)$ of a group $G$ on a vector space $V$, we will use the word {\em representation} to refer freely to either $\rho$ or $V$.

\begin{definition}\label{def:pnc}
If $\rho:G\rightarrow GL(V)$ is a representation of $G$ and $x,y\in G$ are such that $\rho(x)$ and $\rho(y)$ do not share any eigenvectors, then we say $x,y$ \textbf{noncommute purely in $V$} (or $\rho$), and $V$ (or $\rho$) is \textbf{purely noncommuting}, or \textbf{PNC}, \textbf{for} the pair $x,y$. We also say $G$ itself is \textbf{PNC for} $x,y$ if it has a representation that is PNC for $x,y$.

Finally, a group $G$ is \textbf{purely noncommuting}, or \textbf{PNC}, if it is PNC for every noncommuting pair.
\end{definition}

We may also be interested to know if we can find a single representation that has this property for all of $G$'s noncommuting pairs. Therefore we make a second definition:

\begin{definition}\label{def:PNC}
If there exists a representation $\rho:G\rightarrow GL(V)$ of $G$ such that all noncommuting pairs $x,y\in G$ noncommute purely in $\rho$, then we say $G$ is \textbf{strongly purely noncommuting} or \textbf{SPNC}.
\end{definition}

%In this language, a group $G$ is PNC if for each noncommuting pair one can find a PNC representation, and it is SPNC if one can find a single representation that is PNC for all pairs.

There is an additional, geometric motivation for these notions. When one takes the quotient of a smooth complex algebraic variety $X$ by the action of a finite group $G$, the resulting variety typically has singular points. By Hironaka's theorem, the singularities can be resolved by a sequence of blowups. However, in general, it is a hard problem to make the desingularization process completely constructive.

On the other hand, if the singularities are {\em abelian}, meaning that they are locally isomorphic to the quotient of $\C^n$ by a finite abelian group, then the desingularization can be accomplished in an explicit way (see \cite{cox}, Chapters 10 and 11). Thus abelian singularities are mild from the point of view of resolution of singularities.

%\begin{definition}
%If $V$ is an algebraic variety over $\C$, $v\in V$ is an \textbf{abelian quotient singularity} (or simply \textbf{abelian singularity}) if the completion 
%\[
%\widehat{\mathcal{O}_{V,v}}
%\]
%of the local ring at $v$ is isomorphic to the completion 
%\[
%\widehat{\mathcal{O}_{Y,0}},
%\]
%where $Y$ is the quotient of $\C^n$ by a finite abelian group $A$, acting linearly.
%\end{definition}

The singular points of the quotient $X/G$ are automatically abelian if for any point $p\in X$, the point stabilizer $G_p$ is abelian. In this case, the image of $p$ in the quotient $X/G$ is locally isomorphic to the quotient of $\C^n$ by $G_p$. (This is not a necessary condition for the quotient singularities to be abelian; see section \ref{sec:absing} below.)

Let $G$ be a finite group with a representation $V$. Then $G$ acts on the projective space $\Proj(V)$. The SPNC property guarantees that $\Proj(V)/G$ will have at worst abelian singularities, as follows. If any two $x,y\in G$ both stabilize $p\in \Proj(V)$, then a representative of $p$ in $V$ is precisely a common eigenvector for the actions of $x$ and $y$ on $V$. Thus if $V$ realizes $G$ as SPNC, it must be that $x,y$ commute. It follows that $G_p$ is abelian.

PNC groups themselves have a related property. If $G$ is a PNC group, and $X = \prod\Proj(V_i)$, where the product is taken over the irreducible representations of $G$, then any point stabilizers for the action of $G$ on $X$ are abelian, by similar reasoning, and therefore any singular points of $X/G$ are abelian.

Thus SPNC and PNC groups have actions, respectively, on projective spaces and products thereof, whose quotients have mild singularities. %When such a group has a center, the action is not faithful, but it can be made faithful by replacing the variety carrying the action with a suitable blowup.

This article is an investigation into PNC and SPNC groups. 
%Section \ref{sec:example} works out a concrete example. 
Section \ref{sec:preliminaries} sets up notation and conventions, and makes some first observations used in the sequel. In section \ref{sec:supersolvable}, we show that supersolvable groups are always PNC (theorem \ref{thm:supersolvable}). At the other end of the spectrum, we show in section \ref{sec:nonabsimple} that nonabelian simple groups are never PNC (theorem \ref{nonabsimple}). In section \ref{sec:SPNCfam}, we prove that a group is SPNC if all its nonabelian subgroups are sufficiently large in a suitable sense (proposition \ref{bigirrep}). Section \ref{sec:metabelian} investigates the PNC property for metabelian groups. We show that a metabelian group is necessarily PNC if its commutator subgroup's cyclic prime-power-order factors are all different (theorem \ref{thm:metabcrit}), and also exhibit an infinite family of metabelian groups that are not PNC (theorem \ref{affinenotPNC}). In section \ref{sec:permutations}, we give a group-theoretic characterization of when two permutations noncommute purely in the standard representation of the symmetric group (proposition \ref{prop:permutations}). We close with further questions in section \ref{sec:questions}.

\section{Notation and preliminaries}\label{sec:preliminaries}

Throughout, for commutators and conjugates we adopt the right-action notation $x^y = y^{-1}xy$ and $[x,y]=x^{-1}y^{-1}xy = x^{-1}x^y$.

%We use the word {\em character} in two well-established but distinct senses. When speaking of a representation $\rho:G\rightarrow GL(V)$ of arbitrary degree of a {\em nonabelian} group $G$, \textbf{the character of the representation} is the class function $\chi:G\rightarrow\C$ giving traces of the actions of elements of $G$ on the representation space. In other words, 
%\begin{align*}
%\chi:G &\rightarrow \C \\
%g &\mapsto \Tr \rho(g).
%\end{align*}
%We will also refer to $\chi$ in this situation as \textbf{a character of the group $G$}. If $\rho$ is an irreducible representation of $G$, then $\chi$ is called an \textbf{irreducible character} of $G$.

%When speaking of an {\em abelian} group, a \textbf{character of the group} is an element of its Pontryagin dual, i.e. a one-dimensional representation of the group.

%\begin{remark}
%A one-dimensional representation 
%\[
%\chi: A \rightarrow GL(1,\C)
%\]
%of an abelian group $A$ is a character in both of the above senses. By identifying $GL(1,\C)$ with $\C^\times\subset \C$, it is a class function giving its own traces; but it is also an element of the pontryagin dual, which emphasizes that it has the added feature of being multiplicative: for $x,y\in A$, we have
%\[
%\chi(xy) = \chi(x)\chi(y).
%\]
%\end{remark}

%The word {\em character} is sometimes also used to refer to one-dimensional representations of groups that are not abelian. We will avoid this usage entirely and merely call them ``one-dimensional representations."

\begin{notation}\label{not:G/A}
If $G$ is a group and $A\subseteq G$ is a (not necessarily normal) subgroup, we use the symbol $G/A$, as usual, to mean the left coset space  of $A$ in $G$. Then the statement
\[
[s]\in G/A
\]
should be interpreted to mean that $[s]$ is a coset and $s$ is some representative in $G$ of this coset. We will then sometimes write
\[
g^s
\]
to mean the conjugate of $g$ by {\em any representative} of the coset $[s]$. We will only use this notation when the setting renders the choice of coset representative inconsequential. The main examples are below in \ref{rmk:extendedchipreservesmult} and \eqref{eq:resofind}.

%Note that when $A$ is normal, this notation is consistent with the meaning of $G/A$ as the quotient group.
\end{notation}

\begin{notation}\label{not:extendedchi}
If $\chi$ is a character of a group $A$ that is embedded in a larger group $G$, we adopt the convention that $\chi$ can be extended to a function on $G$, also called $\chi$, by assigning it the value $0$ outside $A$. More precisely, define a new function $\overline{\chi}$ by
\[
\overline{\chi}(g) =
\begin{cases} 
\chi(g), & g\in A\\
0,& g\notin A,
\end{cases}
\]
and then set $\chi = \overline{\chi}$. Note that $\chi$ is not necessarily a character of $G$.
\end{notation}

\begin{observation}\label{rmk:extendedchipreservesmult}
Convention \ref{not:extendedchi} allows us to write the formula 
\[
\Ind_A^G\chi(g) = \sum_{[s]\in G/A} \chi(g^s),
\]
giving the character of an induced representation. Per \ref{not:G/A}, this formula does not depend on the choice of coset representative $s\in[s]$: if $g^s\notin A$, then 
\[
g^{sa} = \left(g^s\right)^a\notin A
\]
either, so $\chi(g^s) = \chi(g^{sa})=0$, while if $g^s\in A$, then 
\[
\chi(g^{sa}) = \chi \left(\left(g^s\right)^a\right) = \chi(g^s)
\]
because $\chi$ is a class function on $A$.

In the case that $A$ is abelian and $\chi$ is multiplicative on $A$, the extended meaning of $\chi$ given by \ref{not:extendedchi} preserves the multiplicativity relation $\chi(gh)=\chi(g)\chi(h)$ as long as at least one of $g,h$ is in $A$. For if one of $g,h$ is in $A$ while the other is not, then $gh$ is not in $A$, so that $\chi(gh) = 0 = \chi(g)\chi(h)$.
\end{observation}

If $G$ fails to be PNC, then it means that there is a noncommuting pair $x,y\in G$ such that in every representation of $G$, $x$ and $y$ share a common eigenvector. This is equivalent to the statement that any representation of $G$, when restricted to the nonabelian subgroup $H$ generated by $x$ and $y$, will contain some one-dimensional representation of $H$. This is a fact about $H$ that does not depend on the choice $x,y$ of its generators. Any other $x',y'$ that also generate $H$ will also obstruct PNCness, i.e. they will share a common eigenvector in every representation of $G$. 

Conversely, if $G$ is PNC, then for every pair $x,y$ of noncommuting elements, there is a representation $V$ in which they do not share a common eigenvector. This means that the restriction of $V$ to $H=\langle x,y\rangle$ must not contain any one-dimensional representations of $H$. Again, this is a statement about $V$ and $H$ that does not depend on the choice of generators $x,y$ for $H$.

These considerations motivate the following definition:

\begin{definition}\label{def:sepforH}
Let $H\subset G$ be a nonabelian subgroup. Given a representation $\rho$ (respectively $V$) of $G$, we say $\rho$ (respectively $V$) is \textbf{PNC for} $H$ if $\rho$'s (respectively $V$'s) restriction to $H$ does not contain any one-dimensional representations of $H$. We also say $G$ is \textbf{PNC for $H$} if $G$ has a representation $V$ that is PNC for $H$.
\end{definition}

It is not hard to see that a group $G$ is PNC if and only if it is PNC for each of its 2-generated nonabelian subgroups. In fact the quantification in this statement can be restricted: $G$ is PNC if and only if it is PNC for each of its {\em minimal nonabelian subgroups}, i.e. its nonabelian subgroups whose proper subgroups are all abelian. (Such groups are automatically 2-generated.)

Among the minimal nonabelian groups are the dihedral groups of order $2p$, $p$ an odd prime, and $8$. Dihedral groups play an important role in a number of our arguments because they are particularly adept at obstructing PNCness, by the following lemma and corollary \ref{D4obstruct}. Per the dihedral group's interpretation as symmetries of a regular polygon, we refer to its cyclic, index-2 subgroup as the \textbf{rotation subgroup} and the elements of the nontrivial coset of this subgroup as \textbf{reflections}.

%, so we take a moment to recall the definition and highlight the property that will be useful to us.

%\begin{definition}\label{def:dihedral}
%The \textbf{dihedral group} $D_n$ of order $2n$ is the nontrivial semidirect product 
%\[
%C_n\rtimes C_2 = \langle r,f\mid r^n=f^2= 1, r^f = r^{-1}\rangle.
%\]
%Note that we must have $n\geq 3$ or else the action of $f$ on $r$ is trivial. The \textbf{rotation subgroup} is the subgroup $C_n = \langle r\rangle$, and the \textbf{reflections} are the elements of the rotation subgroup's nontrivial coset $\langle r\rangle f$. 
%\end{definition}

%\begin{remark}
%These names come from the description of $D_n$ as the group of Euclidean symmetries of a regular $n$-gon ({\em dihedron}) in $\R^2$. One may take $r$ to be the $2\pi/n$-rotation of the $n$-gon about its center, and $f$ to be any reflection symmetry.
%\end{remark}

%Recall that a representation that splits into multiple copies of the same irreducible representation is said to be \textbf{isotypical}.

\begin{lemma}\label{lem:dihedral}
Suppose $G$ contains a dihedral group $D$, and $G$ has no irreducible representation whose character simultaneously (i) is identically zero on $D$'s reflections, and (ii) sums to zero on $D$'s rotation subgroup. Then $G$ is not PNC.
\end{lemma}

\begin{proof}
Every irreducible representation of $D$ of degree greater than 1 has properties (i) and (ii). In fact, they are all of degree 2, and (i) and (ii) are immediate consequences of the fact (\cite[pp.~35--37]{serre}) that they are all induced from nontrivial one-dimensional representations of the rotation subgroup. Since these properties are both linear, they also hold for any representation of $D$ that does not contain a one-dimensional representation. In particular, the hypothesis guarantees that any irreducible representation of $G$ (and therefore any representation at all) will contain a one-dimensional subrepresentation when restricted to $D$. Thus, $G$ is not PNC for $D$, and therefore not PNC.
\end{proof}

We use this lemma in several proofs in sections \ref{sec:nonabsimple} and \ref{sec:metabelian}.

\section{Supersolvable groups}\label{sec:supersolvable}

A finite solvable group admits a normal series with abelian quotients (the derived series), and a subnormal series with cyclic quotients (any composition series). If we strengthen this requirement to a normal series with cyclic quotients, we can guarantee that the group is PNC. Recall that groups with such a normal series are called \textbf{supersolvable}.

A well-known formula (\cite[Proposition 22]{serre}) gives a decomposition of the restriction to a subgroup of a representation induced from another subgroup. In the case that the two subgroups coincide and are normal, 
%say $N$ is a normal subgroup of a group $G$, and $\rho: N\rightarrow GL(V)$ is a representation of $N$, 
then this formula takes the form
\begin{equation}\label{eq:resofind}
\Res^G_N\Ind_N^G \rho \cong \bigoplus_{[s]\in G/N} \rho^s,
\end{equation}
where each $\rho^s:N\rightarrow GL(V)$ is defined by $\rho^s(x) = \rho(x^s)$.

%\begin{remark}
%The formula \eqref{eq:resofind} illustrates the point discussed in \ref{not:G/A}. As a map, $\rho^s: N\rightarrow GL(V)$ may depend on the choice of representative $s$ for a coset $[s]\in G/N$. However, if $s' = sn$ with $n\in N$, we have
%\[
%\rho^{s'}(x) = \rho(x^{sn}) = \rho(n^{-1}x^sn) = \rho(n^{-1})\rho(x^s)\rho(n) = \rho(n)^{-1}\rho^s(x)\rho(n),
%\]
%so that $\rho^{s'}$ and $\rho^s$ are isomorphic representations of $N$, with the isomorphism induced by the automorphism $\rho(n)\in GL(V)$. Thus the right side of \eqref{eq:resofind} is unambiguous up to isomorphism of representations, which is what the context requires.
%\end{remark}

\begin{observation}\label{lem:cyclic1dfaithful}
Cyclic groups are precisely those finite groups possessing a faithful one-dimensional representation, since every finite subgroup of $\C^\times$ is cyclic.
\end{observation}

%\begin{proof}
%For a cyclic group of order $n$, one obtains a one-dimensional faithful representation by mapping a generator to an $n$th root of unity $\zeta_n$.

%In the other direction, the image of any one-dimensional representation of a finite group $G$ is cyclic, since all finite subgroups of $\C^\times$ are cyclic. If the representation is faithful, this means $G$ is cyclic.
%\end{proof}

\begin{observation}\label{lem:commutator1}
If $V$ is a vector space and $x,y\in GL(V)$ share an eigenvector $v$, then their commutator $[x,y]\in GL(V)$ also shares this eigenvector, and it has eigenvalue $1$, i.e.
\[
v[x,y] = v.
\]
\end{observation}

%\begin{proof}
%Suppose the eigenvalues of $X,Y$ corresponding to the eigenvector $v$ are $\alpha,\beta$. We have
%\[
%[X,Y]v = X^{-1}Y^{-1}XYv = \alpha^{-1}\beta^{-1}\alpha\beta v = v.\qedhere
%\]
%\end{proof}

\begin{thm}\label{thm:supersolvable}
Let $G$ be a finite group admitting a normal series with cyclic quotients, i.e. a finite supersolvable group. Then $G$ is PNC.
\end{thm}

\begin{proof}[Proof of Theorem \ref{thm:supersolvable}]
Let $x,y$ be any two noncommuting elements. Then $[x,y]$ is a nontrivial element of $G$. Let
\[
G=G_0\triangleright G_1\triangleright\dots\triangleright G_n = \{1\}
\]
be the assumed normal series with cyclic quotients. Since $[x,y]$ is nontrivial, it is in $G_i\setminus G_{i+1}$ for some $i=0,\dots,n-1$. By assumption, $G_i/G_{i+1}$ is cyclic, so there exists a faithful one-dimensional representation 
\[
\rho : G_i/G_{i+1}\rightarrow \C^\times
\]
by observation \ref{lem:cyclic1dfaithful}. Precomposing with the canonical homomorphism 
\[
\pi: G_i\rightarrow G_{i}/G_{i+1},
\]
we obtain a one-dimensional representation $\phi =\rho \pi$ of $G_i$ that is nontrivial outside of $G_{i+1}$.

We will now show that the induced representation
\[
\Phi = \Ind_{G_i}^G\phi
\]
is PNC for $x,y$. This will be done by showing that $\Phi([x,y])$ does not have $1$ as an eigenvalue. It will then follow that $\Phi(x),\Phi(y)$ do not share an eigenvector, for if they did, their commutator $\Phi([x,y])$ would have $1$ as an eigenvalue, by observation \ref{lem:commutator1}. As $x,y$ are an arbitrary noncommuting pair, this will complete the proof that $G$ is PNC.

Since $G_i$ is normal in $G$, \eqref{eq:resofind} tells us that
\[
\Phi|_{G_i} = \bigoplus_{[s]\in G/G_i} \phi^s.
\]
Since $\phi$ is one-dimensional, $\phi^s$ is as well, for each $s$, so that this formula splits $\Phi$ into one-dimensional representations on $G_i$. It follows that for any given $g\in G_i$, the eigenvalues of $\Phi(g)$ are just the values of the $\phi^s(g)\in \C^\times$.

We apply this with $g=[x,y]\in G_i$. By assumption, $[x,y]$ lies outside of $G_{i+1}$. Since $G_{i+1}$ is normal in $G$, $[x,y]^s$ also lies outside of $G_{i+1}$ for each $s$. Since $\phi$ is nontrivial outside of $G_{i+1}$ by construction, this means that $\phi^s([x,y]) = \phi([x,y]^s)$ is not equal to $1$ for any $s$. Thus no eigenvalue of $\Phi([x,y])$ is $1$. This completes the proof.
\end{proof}

\begin{cor}
Finite nilpotent groups are PNC.
\end{cor}

\begin{proof}They are supersolvable.\end{proof}

Recall that a group $G$ is called \textbf{cyclic-by-abelian} if it has a cyclic normal subgroup $C$ such that the quotient $G/C$ is abelian.

\begin{cor}\label{cyclicbyabelian}
If a finite group $G$ is cyclic-by-abelian, then it is SPNC.
\end{cor}

\begin{proof}
Let $C\triangleleft G$ be cyclic with $G/C$ abelian. Then $[G,G]\subset C$, thus every nontrivial commutator is in $C\setminus\{1\}$. There exists a character $\phi$ of $C$ that is nontrivial on $C\setminus\{1\}$, by observation \ref{lem:cyclic1dfaithful}. Then $\Phi = \Ind_{C}^G \phi$ is a representation in which, by the exact same argument as in the proof of theorem \ref{thm:supersolvable}, no nontrivial commutator has 1 as an eigenvalue, and therefore in which no pair of noncommuting elements shares an eigenvector. Thus $\Phi$ is PNC for any noncommuting $x,y$, so it manifests $G$ as SPNC.
\end{proof}

There is no hope of a similar result about groups which are merely solvable:

\begin{prop}\label{prop:S4}
The symmetric group $S_4$ is not PNC.
\end{prop}

This can be proven by direct reference to $S_4$'s character table, 
%as was done for $A_5$ in section \ref{sec:example}, 
but the following more conceptual proof contains a device that will be used again later.

\begin{proof}
Consider the subgroup $D_4 = \langle (1234),(13)\rangle \subset S_4$. We will show that no representation of $S_4$ is PNC for this subgroup.

$D_4$ acts faithfully on the plane as the symmetry group of a square. This is its only irreducible representation of degree greater than one. 
%(Proof: it has five conjugacy classes, thus five irreducible representations. It has three subgroups of index 2, implying three nontrivial homomorphisms to $\{\pm 1\}$. These, and the trivial representation, account for four of the five.) 
The character $\chi$ of this representation is given by
\[
\begin{array}{ c | c c c c c }
 & 1 & (13)(24) & (1234) & (13) & (12)(34) \\
\hline \chi & 2 & -2 & 0 & 0 & 0
\end{array}
\]
Notice that $\chi$ separates the central element $(13)(24)$ from the class of the reflection $(12)(34)$. On the other hand, in $S_4$ these elements are conjugate. Therefore no class function on $S_4$, in particular no character of $S_4$, can separate them. It follows that no character of $S_4$ restricts to a multiple of $\chi$; thus the restriction to $D_4$ of any representation of $S_4$ must contain some one-dimensional representation of $D_4$, so no representation of $S_4$ is PNC for $D_4$.
\end{proof}

The only feature of $S_4$ used in this proof is that it contains $D_4$ in such a way that the central involution is conjugate to one of the other involutions. Therefore the argument generalizes:

\begin{cor}[$D_4$ obstruction]\label{D4obstruct}
If a finite group $G$ contains $D_4$ in such a way that the nontrivial central element in $D_4$ is conjugate to one of the other involutions, then $G$ is not PNC.\qed
\end{cor}

\section{Nonabelian simple groups}\label{sec:nonabsimple}

While section \ref{sec:supersolvable} shows that there are plenty of PNC groups, there are also plenty of groups which are not PNC. Recall that a family of objects $\mathscr{F}$ is said to be \textbf{upward-closed} if whenever an object $A$ is in $\mathscr{F}$ and embeds in an object $B$, then $B$ is in $\mathscr{F}$ too.

\begin{lemma}\label{lem:upwardclosed}
The family of non-PNC groups is upward-closed.
\end{lemma}

\begin{proof} 
If a noncommuting pair in a group $G$ shares an eigenspace in every representation of $G$, it also does so in every representation $\rho$ of any group containing $G$, since $\rho$ is also a representation of $G$ by restriction. Thus if $G$ is not PNC for $x,y\in G$, no overgroup of $G$ can be PNC for $x,y$ either.
\end{proof}

Thus by proposition %s \ref{prop:A5} and 
\ref{prop:S4}, no group containing $S_4$ % or $A_5$
is PNC.  But more broadly:

\begin{thm}\label{nonabsimple}
No nonabelian simple group is PNC.
\end{thm}

This theorem is the main goal of this section. The structure of the proof is as follows. By a 1997 result of Barry and Ward (\cite[Theorem 1]{barryward}), every nonabelian simple group contains a \textbf{minimal simple group}, i.e. a nonabelian finite simple group all of whose proper subgroups are solvable. Such groups were classified in 1968 by Thompson, and they are all of the form $PSL(2,q)$ for $q$ a prime power $\geq 4$, $Sz(2^p)$ for $p$ an odd prime, or $PSL(3,3)$ (\cite[Corollary 1]{thompson}). We will show that none of these groups is PNC by giving explicit dihedral subgroups for which they are not PNC. The result for all nonabelian simple groups will then follow by lemma \ref{lem:upwardclosed}. Here are the precise details:

%\begin{definition}
%A \textbf{minimal simple group} is a nonabelian finite simple group all of whose proper subgroups are solvable.
%\end{definition}

%\begin{lemma}[\cite{barryward}, Theorem 1]\label{lem:barryward}
%Every nonabelian finite simple group contains a minimal simple group.
%\end{lemma}

%\begin{lemma}[\cite{thompson}, Corollary 1]\label{lem:thompson}
%Every minimal simple group is among the following:
%\begin{enumerate}
%\item $PSL(2,q)$ for $q$ a prime power $\geq 4$.
%\item The Suzuki group $Sz(2^p)$ for $p$ an odd prime.
%\item $PSL(3,3)$.
%\end{enumerate}
%\end{lemma}

\begin{remark}
The group $PSL(2,q)$ is only minimal simple for certain $q$. The version of the statement in \cite{thompson} is sharper than what we have quoted.
\end{remark}

\begin{lemma}\label{psl2q}
For a prime power $q\geq 4$, $PSL(2,q)$ is not PNC.
\end{lemma}

In fact, we can already know this for $q=\pm 1\mod 8$ since in this case $PSL(2,q)$ contains $S_4$. 
%, and for $q=\pm 1\mod 10$ since in this case it contains $A_5$
One could hope to proceed to the remaining cases. We give a more uniform proof, although some case analysis is inevitable because the representation theory of $PSL(2,q)$ depends on $q$ mod $4$.

\begin{proof}
We will show that the PNC property is obstructed by a dihedral group of order $q-1$ or $q+1$, if $q$ is odd, or $2(q - 1)$, if $q$ is even. It is well-known that $PSL(2,q)$ contains dihedral subgroups of these orders (\cite{dickson}, \S 246, or \cite{king}, Theorem 2.1(d)-(i)). Let $D$ be a dihedral subgroup of $PSL(2,q)$, of order to be specified shortly.

We will show that the order of $D$ can always be chosen so that no irreducible character of $PSL(2,q)$ simultaneously is identically zero on the reflections and sums to zero on the rotation subgroup. Then the desired result will follow from lemma \ref{lem:dihedral}.

The relevant part of the character table of $PSL(2,q)$ is given in table \ref{tbl:psl2q}. The notation below is explained in the caption.

\begin{table}
\[
\text{Case $q=1$ mod $4$.}
\]
\[
\begin{array}{c | c c c }
& 1 & a^\ell & b^m \\
\hline \text{Triv} & 1 & 1 & 1\\
\psi & q & 1 & -1\\
\chi_i & q+1 & \rho^{i\ell}+\rho^{-i\ell} & 0\\
\theta_j& q-1 & 0 & -\left(\sigma^{jm} + \sigma^{-jm}\right)\\
\xi_1 & (q+1)/2 & (-1)^\ell & 0\\
\xi_2 & (q+1)/2 & (-1)^\ell & 0
\end{array}
\]
\vspace{2mm}
\[
\text{Case $q=3$ mod $4$.}
\]
\[
\begin{array}{c | c c c }
 & 1 & a^\ell & b^m\\
\hline \text{Triv} & 1 & 1 & 1\\
\psi & q & 1 & -1\\
\chi_i & q+1 & \rho^{i\ell}+\rho^{-i\ell} & 0\\
\theta_j& q-1 & 0 & -\left(\sigma^{jm} + \sigma^{-jm}\right)\\
\eta_1 & (q+1)/2 & 0 & (-1)^{m+1}\\
\eta_2 & (q+1)/2 & 0 & (-1)^{m+1}
\end{array}
\]
\vspace{2mm}
\[
\text{Case $q$ even.}
\]
\[
\begin{array}{ c | c c c c }
 & 1 & c & a^\ell & b^m\\
\hline \text{Triv} & 1 & 1 & 1 & 1\\
\psi & q & 0 & 1 & -1\\
\chi_i & q+1 & 1 & \rho^{i\ell} + \rho^{-i\ell} & 0\\
\theta_j & q-1 & -1 & 0 & - \left(\sigma^{jm} + \sigma^{-jm}\right)
\end{array}
\]
\caption{Character table of $PSL(2,q)$. The symbols $\rho,\sigma$ are primitive $(q-1)$th and $(q+1)$th roots of unity respectively. For odd $q$, respectively even $q$, $a$ is the class of elements of order $(q-1)/2$, respectively $q-1$, and $b$ is the class of elements of order $(q+1)/2$, respectively $q+1$. For odd $q$, $i$ and $j$ are even integers and we have omitted the classes of elements of order dividing $q$. For even $q$, $i,j$ are integers and $c$ is the class of involutions. Source: \cite{dornhoff}, \S 38.}\label{tbl:psl2q}
\end{table}

There are three cases to consider: $q=1$ mod $4$, $q=3$ mod $4$, and $q$ even.

In the case $q=1\mod 4$, the class of involutions is $a^{(q-1)/4}$ in the table. Take $D$ to be of order $q+1$, so its rotation subgroup, of order $(q+1)/2$, consists of the identity and elements in the classes $b^m$. The only irreducible characters of $PSL(2,q)$ that are zero on $a^{(q-1)/4}$ are those of the cuspidal representations $\theta_j$. Their absolute value is $q+1$ on the identity and is bounded by $2$ on the elements $b^m$. Thus the sum of any of these characters across the rotation subgroup has absolute value bounded below by
\[
q+1 - 2\left(\frac{q+1}{2} - 1\right) =2>0.
\]
Therefore no irreducible character of $PSL(2,q)$ is simultaneously zero on $D$'s reflections and sums to zero on $D$'s rotation subgroup.

The other two cases are similar. 

For $q=3\mod 4$, one uses $D$ of order $q-1$ instead of $q+1$. This time, the nontrivial rotations are in the classes $a^\ell$, and the class of involutions is $b^{(q+1)/4}$. Only the principal series characters $\chi_i$ are zero on this latter class. Practically the same calculation shows their sum across the rotation subgroup is bounded below by $2$.

%In the case $q=3\mod 4$, the class of involutions is $b^{(q+1)/4}$. Take $D$ to be of order $q-1$, so its rotation subgroup consists of the identity and elements $a^\ell$. In this case, only the principal series characters $\chi_i$ are zero on the class of involutions, and they are $q-1$ on the identity and bounded by $2$ in absolute value on the nontrivial rotations. Thus these characters' sums across $D$'s rotation subgroup has absolute value again bounded below by
%\[
%q-1 - 2\left(\frac{q-1}{2} - 1\right) = 2>0
%\]
%and so cannot be PNC for $D$.

In the final case of even $q$, take $D$ to be of order $2(q-1)$, so the nontrivial rotations are $a^\ell$ and the involutions are $c$. The only irreducible character that is zero on $c$ is the Steinberg character $\psi$, which is positive on the $a^\ell$'s, so the sum across the rotation subgroup of $D$ is positive.
\end{proof}

%\begin{remark}
%This proof generalizes the proof for $A_5$ given above in proposition \ref{prop:A5} in two different ways. In that proof, the PNC-obstructing subgroup was $S_3$, which is isomorphic to the dihedral group of order $6$. Noting that $A_5\cong PSL(2,4)$, the even case, this proof selects the dihedral group of order $2(q-1) = 2\cdot 3 = 6$ to obstruct PNCness. Alternatively, noting that $A_5\cong PSL(2,5)$, the $1$ mod $4$ case, it selects the dihedral group of order $q+1 = 5+1 = 6$ as well.
%\end{remark}

\begin{lemma}\label{psl33}
$PSL(3,3)$ is not PNC.
\end{lemma}

\begin{proof}
The argument is identical to that given for $S_4$ (proposition \ref{prop:S4}). $PSL(3,3)=SL(3,3)$ contains a subgroup $D_4$ generated by
\[
r = \begin{pmatrix} &-1& \\1& & \\ & &1\end{pmatrix},\;f = \begin{pmatrix}1& & \\ &-1& \\ & &-1\end{pmatrix}
\]
The central involution in this copy of $D_4$ is
\[
\begin{pmatrix}-1& & \\ &-1& \\ & &1\end{pmatrix}
\]
But this is conjugate in $PSL(3,3)$ to $f$, so apply corollary \ref{D4obstruct}.
\end{proof}

\begin{lemma}\label{suzuki}
The Suzuki group $Sz(2^p)$ ($p$ an odd prime) is not PNC.
\end{lemma}

\begin{proof}
The proof is essentially the same as that for $PSL(2,q)$ when $q$ is even.

Like $PSL(2,q)$, the Suzuki group $Sz(q), \;q = 2^p$ has a single conjugacy class of involutions (\cite[Proposition 7]{suzuki}). It also contains a dihedral group $D$ of order $2(q-1)$: when $Sz(q)$ is realized as a permutation group as in Suzuki's original presentation, this is the normalizer of the stabilizer of two points (\cite[Proposition 3]{suzuki}).

The character table of $Sz(q)$ is given in table \ref{tbl:Szq}, which is taken from \cite[Theorem 13]{suzuki}.

\begin{table}
\[
\begin{array}{ c | c c c c c c }
\text{Class name:} & 1 & \sigma & \rho,\rho^{-1} & \pi_0 & \pi_1 & \pi_2\\
\text{Order divides:} & 1 & 2 & 4 & q-1 & q + r + 1 & q - r + 1\\
\hline X & q^2 & 0 & 0 & 1 & -1 & -1 \\
X_i & q^2 + 1 & 1 & 1 & \varepsilon_0^i(\pi_0) & 0 & 0\\
Y_j & (q - r + 1)(q-1) & r-1 & -1 & 0 & -\varepsilon_1^j(\pi_1) & 0\\
Z_k & (q + r + 1)(q-1) & -r-1 & -1 & 0 & 0 & -\varepsilon_2^k(\pi_2)\\
W_l & r(q-1)/2 & -r/2 & \pm ri/2 & 0 & 1 & -1
\end{array}
\]
\caption{Character table of $Sz(q)$. Here, $r = \sqrt{2q}$. The classes called $\pi_0,\pi_1,\pi_2$ consist of elements belonging to certain cyclic subgroups $A_0,A_1,A_2$, and the $\varepsilon$'s are certain characters of these subgroups.}\label{tbl:Szq}
\end{table}

%Per lemma \ref{dihedralreps}, as in the proof for $PSL(2,q)$, in order for an irreducible representation of $Sz(q)$ to be PNC for $D$ its character would have to be zero on the class of involutions and sum to zero on the elements of $D$'s rotation subgroup. 
There is only one irreducible character of $Sz(q)$ that is zero on the class of involutions ($X$ in the table), and it is positive on all the elements of $D$'s rotation subgroup (which, besides the trivial class, are in the classes Suzuki calls $\pi_0$, as their orders divide $q-1$). So, by lemma \ref{lem:dihedral}, $Sz(q)$ is not PNC. 
\end{proof}

\begin{proof}[Proof of theorem \ref{nonabsimple}]
As noted above, \cite[Theorem 1]{barryward} and \cite[Corollary 1]{thompson} give us that every nonabelian finite simple group contains $PSL(2,q)$ for $q\geq 4$, $Sz(2^p)$ for $p$ odd prime, or $PSL(3,3)$. None of these is PNC by \ref{psl2q}, \ref{psl33}, and \ref{suzuki}, so \ref{lem:upwardclosed} then implies that no nonabelian finite simple group is PNC.
\end{proof}

\section{A family of SPNC groups}\label{sec:SPNCfam}

The argument of theorem \ref{thm:supersolvable} shows supersolvable groups are PNC by finding, for any noncommuting pair, a representation in which its commutator does not have 1 as an eigenvalue. A group can be PNC without this. For example, $A_4$ is a group in which, in {\em every} representation, {\em every} commutator has $1$ as an eigenvalue. Nonetheless, the standard 3-dimensional representation of $A_4$ actually realizes it as SPNC: $A_4$ is a minimal nonabelian group, thus any pair of noncommuting elements generates the whole group,
% (\ref{lem:minimal2generated})
and therefore cannot have a common eigenspace in this representation because it is irreducible.

In a similar way one sees immediately that any minimal nonabelian group is SPNC. This is actually a special case of a more general phenomenon that forces a group to be not just PNC but SPNC:

\begin{prop}\label{bigirrep}
Let $G$ be a finite group with an irreducible representation $V$ of degree $d$ that exceeds the index of any of its nonabelian subgroups. Then $V$ realizes $G$ as SPNC.
\end{prop}

\begin{proof}
Let $x,y$ be a pair of noncommuting elements of $G$, and let $H=\langle x,y\rangle$. By assumption, $[G:H]<d$. Now let $L$ be any one-dimensional representation of $H$. Then, by Frobenius reciprocity, the number of times that $L$ occurs in the restriction of $V$ to $H$ is equal to the number of times $V$ occurs in the induced representation $\Ind_H^G L$. Since the dimension of this representation is 
\[
[G:H]<d = \dim V,
\]
this number is zero. So no one-dimensional representation $L$ occurs in the restriction of $V$ to $H$, i.e. $x,y$ do not have a common eigenspace in $V$.
\end{proof}

\begin{remark}
As we have seen, failure to be PNC is always caused by specific obstructing subgroups. For example, the obstruction for $S_4$ is $D_4$ (\ref{prop:S4}). Proposition \ref{bigirrep} shows that a subgroup obstructing PNCness cannot be ``too big." Indeed, $D_4\subset S_4$ is ``as big as possible," since $S_4$ has an irreducible representation (in fact, two) of degree $3$, equal to the index.
\end{remark}

Beyond the trivial case of minimal nonabelian groups, an example of a family of groups satisfying the hypothesis of \ref{bigirrep} is the semidirect product of the additive group $A$ of the field $\F_{p^2}$, for $p$ a prime congruent to $1$ mod $4$, by an automorphism given by multiplication by a field element of multiplicative order $d=(p+1)/2$. This group has an irreducible representation of degree $d$, but every nonabelian subgroup properly contains $A$ and so has index less than $d$. (See \cite[Proposition 1.1.49]{blumsmith} for details.)

\section{Metabelian groups}\label{sec:metabelian}

Recall that a group is called \textbf{metabelian} if it has an abelian normal subgroup with an abelian quotient, in other words if it is solvable of height two. In this section we investigate the PNC property for metabelian groups.

The results of sections \ref{sec:supersolvable}, \ref{sec:nonabsimple}, and \ref{sec:SPNCfam} show that PNCness is loosely correlated with abelianness -- the ``extremely nonabelian" simple groups are never PNC, while ``almost abelianness" of various kinds (nilpotence and supersolvability, per section \ref{sec:supersolvable}, and having all nonabelian subgroups ``large," per section \ref{sec:SPNCfam}) guarantee PNCness. (Abelian groups themselves are PNC, vacuously. It is perhaps a failing of our nomenclature that abelian groups are called ``purely noncommuting".) Based on this intuition, the authors thought that metabelian groups might be always PNC; but this turns out not to be the case.

Recall that, for a prime power $q$, the \textbf{affine group}, or \textbf{affine linear group}, $AGL(n,q)$, is the semidirect product of the additive group of the $\F_q$-vector space $\F_q^n$ by the group $GL(n,q)$ of linear automorphisms of this space.
%group of transformations of the $\F_q$-vector space $\F_q^n$ generated by the linear transformations $GL(n,q)$ and the group $T_{n,q}\cong \F_q^n$ of translations 
%\begin{align*}
%t_x:\F_q^n&\rightarrow \F_q^n\\
%v &\mapsto x+v.
%\end{align*}
%Since conjugation by a linear transformation sends translations to translations, $T_{n,q}$ is a normal subgroup of this group, and since $GL(n,q)$ stabilizes the origin, while $T_{n,q}\cong\F_q^n$ acts freely on the points of $\F_q^n$, their intersection is trivial. Thus
%\[
%AGL(n,q) = \F_q^n \rtimes GL(n,q)
%\]
%is a Frobenius group, with Frobenius kernel $\F_q^n$ and Frobenius complement $GL(n,q)$. In this context one sees $\F_q^n$ as the affine space $\A_{\F_q}^n$, hence the name, because the presence of the translations means one cannot distinguish the origin among the points of $\F_q^n$ from the abstract group action on these points. A choice of origin is equivalent to the choice of a section $GL(n,q)\rightarrow ALG(n,q)$.
For $n=1$, $AGL(1,q) = \F_q^+\rtimes \F_q^\times$ is metabelian. However:

\begin{thm}\label{affinenotPNC}
If $q=p^k$ is an odd prime power with $k>1$, the affine group 
\[
AGL(1,q) = \F_q^+\rtimes \F_q^\times
\]
is not PNC.
\end{thm}

We defer the proof to the end of the section.

In spite of this negative result, a large class of metabelian groups is PNC. If $G$ is metabelian then the commutator subgroup $[G,G]$ is abelian. There is a criterion on the structure of $[G,G]$ that lets us conclude PNCness without knowing anything else about $G$:

\begin{thm}\label{thm:metabcrit}
Let $G$ be a metabelian group and let
\[
[G,G] \cong C_{p_1^{e_1}}\times\dots\times C_{p_k^{e_k}}
\]
be the expression of its commutator as a direct product of cyclic factors of prime power order. If the factors are pairwise nonisomorphic, then $G$ is PNC.
\end{thm}

\begin{remark}
This is not a necessary condition. Many metabelian groups not satisfying the hypothesis of theorem \ref{thm:metabcrit} are still PNC, for example the family of SPNC groups described at the end of section \ref{sec:SPNCfam}. But it shows that PNC metabelian groups are easy to come by.
\end{remark}

The organization of this section is motivated by the proof of \ref{thm:metabcrit}. The fundamental tool is the following lemma, which is of independent utility in investigating the PNC property:

\begin{lemma}[Commutator criterion]\label{commutatorcriterion}
Let $G$ be an arbitrary finite group, $H$ a nonabelian subgroup of $G$, and $V$ a representation of $G$ with character $\chi$. Then $V$ is PNC for $H$ if and only if $\chi$ sums to zero on each coset of $[H,H]$ in $H$.
\end{lemma}

\begin{proof}
The representation $V$ is PNC for $H$ if and only if $V|_H$ contains no one-dimensional representations of $H$. By the orthogonality relations, this is the case if and only if $\chi|_H$ is orthogonal to all of $H$'s one-dimensional representations, with respect to the $H$-invariant inner product
\[
\langle \chi_1,\chi_2\rangle_H = \frac{1}{|H|}\sum_{h\in H} \chi_1(h)\overline{\chi}_2(h).
\]
Now $H$'s one-dimensional representations are precisely the pullbacks to $H$ of all of the characters of the abelian group $H/[H,H]$. These characters span the full space of functions on $H/[H,H]$ (as for any abelian group), so their pullbacks to $H$ span the space of all functions on $H$ constant on each coset of $[H,H]$. The orthogonal complement of this space is clearly the space of class functions that sum to zero on each coset of $[H,H]$, and $V$ is PNC for $H$ if and only if $\chi|_H$ lies in this orthogonal complement.
\end{proof}

Now we begin to assemble the proof of theorem \ref{thm:metabcrit}.

If $G$ is metabelian, then its commutator subgroup $[G,G]$ is abelian, and therefore acts trivially on itself by conjugation. It follows that the conjugation action of $G$ on $[G,G]$ makes the latter a $G/[G,G]$-module.

In what follows we fix the notation that $G$ is a finite metabelian group, $A=[G,G]$, and $Q=G/[G,G]$, so $A$, $Q$ are abelian and $A$ is a $Q$-module.

\begin{prop}\label{subgroupcharacterH}
Let $H\subset G$ be a nonabelian subgroup, so that $K=[H,H]$ is a nontrivial subgroup of $A$. Suppose that $A$ has a character $\chi$ whose kernel does not contain any image of $K$ under the $Q$-action. Then the representation $\Ind_A^G\chi$ is PNC for $H$.
\end{prop}

For example, if $A$ is cyclic (so $G$ is cyclic-by-abelian), it has a character $\chi$ that is a faithful representation of $A$ and is thus nontrivial on all nontrivial subgroups. Therefore $\Ind_A^G\chi$ is PNC for all nonabelian subgroups; thus $G$ is SPNC. This reproduces corollary \ref{cyclicbyabelian}, although without the information (obtained in the proof in section \ref{sec:supersolvable}) that no commutator has an eigenvalue $1$ in this representation.

\begin{proof}[Proof of proposition \ref{subgroupcharacterH}]
We want to show $\Ind_A^G\chi$ is PNC for $H$, and by the commutator criterion (lemma \ref{commutatorcriterion}), this is equivalent to showing that 
\[
\sum_{k\in K}\Ind_A^G\chi(kh) = 0
\]
for all $h\in H$. Actually we even have $\sum_{k\in K}\Ind_A^G \chi(kg)=0$ for every $g\in G$. We see this as follows:
\begin{align*}
\sum_{k\in K}\Ind_A^G \chi(kg) &= \sum_{k\in K}\sum_{[s]\in G/A} \chi((kg)^s)\\
&=\sum_{k\in K} \sum_{[s]\in G/A}\chi(k^sg^s)\\
&=\sum_{k\in K} \sum_{[s]\in G/A}\chi(k^s)\chi(g^s).
\end{align*}
The last equality is because $\chi$, being a character of $A$, is multiplicative, since $k^s\in A$ as $k\in K\subset A$ and $A$ is normal (see \ref{rmk:extendedchipreservesmult}). Reversing the summations and then reindexing the inner sum, we have
\[
\sum_{[s]\in G/A}\left(\sum_{k\in K} \chi(k^s)\right)\chi(g^s) = \sum_{[s]\in G/A}\left(\sum_{k\in K^s} \chi(k)\right)\chi(g^s).
\]
But the inner sum is zero, because by assumption $K^s$ is not contained in $\ker \chi$, therefore $\chi$ restricts to a nontrivial character of $K^s$, and the sum of a nontrivial character over a group is always zero.
\end{proof}

This proposition immediately implies that the following condition on the module structure guarantees PNCness:

\begin{prop}[Subgroup character condition]\label{subgroupcharactercondition}
Suppose $A$ (as $Q$-module) has the property that for any nontrivial subgroup $K\subset A$, $A$ has a character whose kernel does not contain any image of $K$ under the $Q$-action. Then $G$ is PNC.
\end{prop}

\begin{proof}
If $H$ is any nonabelian subgroup of $G$, then proposition \ref{subgroupcharacterH} shows how to construct a representation of $G$ that is PNC for $H$.
\end{proof}

\begin{remark}
In fact, we found theorem \ref{affinenotPNC} by looking for a group where the condition of this proposition fails. 
\end{remark}

The next lemma, whose statement and proof are due to Frieder Ladisch (\cite{ladisch}), links \ref{subgroupcharactercondition} with the hypothesis of theorem \ref{thm:metabcrit}.

\begin{lemma}[Ladisch]\label{ladischlemma}
The following conditions on a finite abelian group $A$ are equivalent:

\begin{enumerate}
\item $A$ has a nontrivial subgroup $K$ such that every character of $A$ is trivial on an image of $K$ under some automorphism of $A$.
\item In a decomposition of $A$ into cyclic factors of prime power order, two of the factors are isomorphic.
\end{enumerate}
\end{lemma}

\begin{proof}
We write $A$ additively.

Condition 2 is fulfilled by $A$ if and only if it is fulfilled by at least one of $A$'s Sylow subgroups. We will show the same for condition 1. If a nontrivial subgroup $K$ of a Sylow subgroup $A_p$ fulfills condition 1 for $A_p$, it also does so for $A$ because automorphisms of $A_p$ extend to $A$ and characters of $A$ restrict to $A_p$. Conversely, if a nontrivial subgroup $K$ of $A$ fulfills condition 1 for $A$, then there is a prime $p$ (any prime dividing $|K|$ in fact) such that $K_p = A_p\cap K$ fulfills condition 1 for $A_p$, since characters of $A_p$ extend to $A$ and automorphisms of $A$ act on $A_p$.

Thus without loss of generality we can suppose $A$ is a $p$-group. If it fulfills condition 2, it has the form $A = F \oplus B$ where $F\cong C_{p^k}\times C_{p^k}$. Then let $K$ be the cyclic subgroup generated by any nonzero element in $F$. Note $\Aut F\subset\Aut A$. Every character of $A$ restricts to a character on $F$ and we assert any character of $F$ is trivial on some $\Aut F$-image of $K$. Indeed, $\Aut F = GL(2,\Z/p^k\Z)$ acts transitively on the elements of $F$ of any given order, and therefore on the order-$|K|$ cyclic subgroups of $F$. Meanwhile every character of $F$ is trivial on some order-$|K|$ cyclic subgroup since it is trivial on some maximal (order $p^k$) cyclic subgroup, as otherwise its image would not be cyclic. This shows (2)$\Rightarrow$(1). 

In the other direction, suppose $K$ fulfills condition 1 for $A$ and consider the subgroup $A_0$ of $A$ of elements of order dividing $p$. This is an $\F_p$-vector space of dimension the $p$-rank of $A$. It has a filtration 
\[
A_0\supset A_1\supset \dots A_k = 0,
\]
where $A_i = A_0\cap p^iA$ is the subgroup of $A_0$ consisting of $p^i$-divisible elements, and $p^k$ is the exponent of $A$. Both $A_0$ and this filtration of it are invariant under automorphisms of $A$. 

Now as $K$ is nontrivial it contains elements of order $p$, so it must meet $A_0$ nontrivially. Thus there is a maximal $i<k$ such that $A_i$ meets $K$ nontrivially; fix this $i$, so that $K\cap A_{i+1} = 0$. Furthermore, as $A_i$ and $A_{i+1}$ are both automorphism invariant, we must have that every image $K'$ of $K$ under $\Aut A$ also meets $A_i$ nontrivially and $A_{i+1}$ trivially.

We assert $A_{i+1}$ is codimension $>1$ in $A_i$. If it were codimension $1$, it would be the kernel of some character on $A_i$, which could be extended to a character $\chi$ of $A$. This character would be nontrivial on every $\Aut A$-image $K'$ of $K$, since they all meet $A_i$ nontrivially outside of $A_{i+1}$. This contradicts the assumption that $K$ fulfills condition 1 for $A$, so we conclude $A_{i+1}$ is codimension $>1$ in $A_i$.

We claim this in turn implies that at least two of $A$'s cyclic factors are isomorphic. Indeed, the dimension of $A_i$ is the number of cyclic factors of $A$ of order at least $p^{i+1}$.\footnote{In fact, if $\lambda = (\lambda_1,\dots,\lambda_r)$ with $\lambda_1\geq \dots\geq \lambda_r$ is the partition describing the type of $A$, so that $A \cong \prod_{\lambda_j} C_{p^{\lambda_j}}$, then the tuple $(\dim A_0,\dots,\dim A_{k-1})$ is the conjugate partition $\lambda'$.} That $A_{i+1}$ is codimension at least two in $A_i$ thus implies that the number of cyclic factors of order at least $p^{i+1}$ is at least two greater than the number of cyclic factors of order at least $p^{i+2}$. This implies that there are at least two cyclic factors of order exactly $p^{i+1}$.

This establishes (1)$\Rightarrow$(2).
\end{proof}

\begin{proof}[Proof of theorem \ref{thm:metabcrit}]
In this situation, by lemma \ref{ladischlemma}, for every subgroup $K$ of $A=[G,G]$, there is a character $\chi$ of $A$ whose kernel does not contain any image of $K$ under $\Aut A$, so $A$ fulfills the hypothesis of proposition \ref{subgroupcharactercondition} for any possible $Q$-action.
\end{proof}

It remains to prove our claim about $AGL(1,q)$. We fix notation:

%\begin{notation}\label{not:AGL}
Let $G = AGL(1,q)$. Then $G = A\rtimes Q$ where $A$ is isomorphic to the additive and $Q$ to the multiplicative group of $\F_q$. The commutator subgroup $[G,G]$ is equal to $A$, so these labels are consistent with those used throughout the section. We identify $G$ as a group of permutations of the  elements of $\F_q$, with $A$ being the translations by $x\in \F_q$ and $Q$ being the multiplications by $a\in \F_q^\times$. Denote the former by $t_x$ and the latter by $m_a$.

\begin{proof}[Proof of theorem \ref{affinenotPNC}]
Recall that $q=p^k$ is an odd, composite prime power.

By \cite[Theorem 6.1]{piatetskishapiro}, the only irreducible representation of $G=AGL(1,q)$ of degree greater than $1$ is the unique representation induced from any nontrivial character of $A$. Its restriction to $A$ is the sum of all nontrivial characters of $A$. Call it $W$. All $G$'s hope of being PNC lies with $W$.

Consider the character $\chi_W$ of $W$. Since $A$ is normal, $\chi_W$ is zero outside of $A$.
% (lemma \ref{inducingfromnormal} of the appendix)
 On $A$, as it is the sum of all nontrivial characters of $A$, it is one less than the sum of {\em all} characters, which is the character of the regular representation. Thus 
\[
\chi_W(1) = |A|-1 = q-1,
\]
and 
\[
\chi_W(t_x) = 0-1 = -1
\]
for all $x\neq 0$ in $\F_q$.

Let $D$ be the subgroup generated by $m_{-1}$ and any $t_x$ with $x\neq 0$. Because $q$ is odd, $m_{-1}$ is the negation map on $\F_q$. Therefore $D$ is a dihedral group of order $2p$, where $p$ is the characteristic of $\F_q$. We now show $W$ is not PNC for $D$.

The rotation subgroup of $D=\langle t_x,m_{-1}\rangle$ is $\langle t_x\rangle$, of order $p$. Thus the sum of $\chi_W$ over the rotation subgroup of $D$ is $1(q-1)+(p-1)(-1) = q-p$. Since $q=p^k$ with $k>1$, this is nonzero, and we can conclude from lemma \ref{lem:dihedral} that $W$ is not PNC for $D$. This concludes the argument.
\end{proof}

\section{Permutations that noncommute purely in the standard representation}\label{sec:permutations}

Although the symmetric group $S_n$ is not PNC for $n\geq 4$ (by \ref{prop:S4} and \ref{lem:upwardclosed}, since it contains $S_4$), we may still be interested, for a given pair of permutations $x,y\in S_n$, whether they noncommute purely in a given representation. Recall that the \textbf{standard representation} $\operatorname{Sta}$ of $S_n$ is the defining permutation representation minus the trivial representation. It can be realized as the action of $S_n$ on the sum-zero subspace $T = \{\sum c_ie_i \mid \sum c_i = 0\}$ of $\C^n$ through permutations of the standard basis vectors $e_i$.

Here we give a group-theoretic characterization of when this representation is PNC for a given $x,y$:

\begin{prop}\label{prop:permutations}
Two permutations $x,y\in S_n$ noncommute purely in the standard representation $\operatorname{Sta}$ of $S_n$ if and only if the subgroup $H=\langle x,y\rangle$ they generate has both of the following properties:
\begin{enumerate}
\item The action of $H$ on the indices $\{1,\dots,n\}$ is transitive.
\item Some point stabilizer $H_j$ (for $j\in \{1,\dots,n\}$) meets every coset of the commutator subgroup $[H,H]$ in $H$.
\end{enumerate}
\end{prop}

\begin{proof}
The defining permutation representation of a permutation group contains one trivial representation for every orbit. Since the standard representation is the defining permutation representation minus one trivial representation, it contains $o - 1$ trivial representations, where $o$ is the number of orbits. Thus $\operatorname{Sta}|_H$ contains no trivial representations if and only if $o=1$, i.e. if and only if $H$ is transitive. Thus failure of (1) implies $x,y$ do not noncommute purely, and it remains to show that (2) is equivalent to noncommuting purely in the presence of (1). 

Thus, we may henceforth assume (1) holds, i.e. $H$ is transitive. In this situation, all point stabilizers $H_j$ are conjugate, and (2) is equivalent to the statement that $H_1$ meets every coset of $[H,H]$. So we need to show that this is equivalent to $\operatorname{Sta}$ being PNC for $H$.

Realize $\operatorname{Sta}|_H$ as $H$'s action on the sum-zero subspace $T\subset \C^n$.

Suppose $v\in T$ lies in a common eigenspace of $x$ and $y$. Then $H$ acts by a one-dimensional representation $\chi:H\rightarrow \C^\times$ on $v$, i.e. $vh = \chi(h)v$ for all $h\in H$.  If
\[
v = \sum c_ie_i,
\]
then this implies $c_1 = c_j\chi(h)$, where $j$ is the image of $1$ under $h$, by matching the $j$th coordinates in $vh = \chi(h)v$. For $h\in H_1$ so that $j=1$, this equation becomes $c_1 = c_1\chi(h)$. It follows that either $c_1=0$, or else $\chi(h)=1$ for all $h\in H_1$. 

Now suppose (2) holds. Then there is an element of $H_1 = \{h\in H\mid e_1h = e_1\}$ in each coset of $[H,H]$. Since $\chi$ is one-dimensional, it factors through the abelianization $H/[H,H]$, thus its value is constant on cosets of $[H,H]$. It follows that either $c_1=0$, or else $\chi(h)=1$ for all $h\in H$. However, both of these conclusions are incompatible with (1) (unless $v=0$). The logic of the first paragraph shows that (1) does not allow $H$ to have a trivial subrepresentation, ruling out $\chi(h)=1,\forall h\in H$ unless $v=0$. Meanwhile, (1) means that for all $j\in\{1,\dots,n\}$, there is an $h\in H$ sending $1$ to $j$, whereupon $c_1 = c_j\chi(h)$ then implies $c_j=0$, so this case is also impossible unless $v=0$. This shows that (1) and (2) preclude a common eigenvector for $x$ and $y$: they noncommute purely.

Conversely, if (2) fails in the presence of (1), then there exists a coset $F$ of $[H,H]$ that does not contain any elements of $H_1$. Since all point stabilizers $H_j$ are conjugate, and since any two conjugate elements lie in the same coset of $[H,H]$, it follows that $F$ consists entirely of fixed-point-free permutations. The value of the character of $\operatorname{Sta}$ is $-1$ on fixed-point-free permutations; thus the sum of $\operatorname{Sta}|_H$ across $F$ is negative. Thus $\operatorname{Sta}$ is not PNC for $H$ by lemma \ref{commutatorcriterion}; i.e. $x,y$ do not noncommute purely in $\operatorname{Sta}$.
\end{proof}

\section{Further questions}\label{sec:questions}

This section collects directions for further inquiry.

\subsection{A classification-free proof for simple groups?}

Our proof of theorem \ref{nonabsimple} rests on the classification theorem for finite simple groups, and thus comes down to case analyses: one, by Barry and Ward (\cite{barryward}), which shows that nonabelian finite simple groups contain minimal simple groups, and another, above in section \ref{sec:nonabsimple}, which shows that minimal simple groups are not PNC. Nonetheless, the case-by-case reasoning used was rather repetitive. This suggests there may be a uniform proof.

\begin{question}
Is there a uniform, classification-free, proof of of theorem \ref{nonabsimple}?
\end{question}

\subsection{Metabelian groups}

The results of section \ref{sec:metabelian} reveal the metabelian case to be richer than we initially expected. As was mentioned in that section, the sufficient criterion \ref{thm:metabcrit} for a metabelian group to be PNC is not close to necessary. A clearer picture of the PNC property for metabelian groups would be desirable.

\begin{question}
Is there a clean if-and-only-if characterization of metabelian groups with the PNC property?
\end{question}

\subsection{Other important representations}

Proposition \ref{prop:permutations} characterizes pairs of purely noncommuting elements group-theoretically in the case of one particularly important representation of one particularly important group. We can ask the same question about other groups and other representations. For example:

\begin{question}
Is there an analogous result to proposition \ref{prop:permutations} for the principal series representations of $PSL(2,q)$?
\end{question}

\subsection{Abelian singularities}\label{sec:absing}

As discussed in the introduction, one motivation for studying the PNC and SPNC properties comes from the search for groups acting on smooth varieties such that the quotients have at worst abelian singularities. But the PNC and SPNC properties are not close to necessary for this. 

First, if a group $G$ has an SPNC central extension $\tilde G$ realized by a representation $V$, then the action of $\tilde G$ on $\Proj(V)$ factors through $G$, and point stabilizers in $G$ are quotients of point stabilizers in $\tilde G$. So in this case $G$ too acts on $\Proj(V)$ with abelian point stabilizers, so $\Proj(V)/G$ has at worst abelian singularities. For example the alternating group $A_5$ is not PNC, but its double cover $\tilde A_5$, the binary icosahedral group, has an SPNC action on $\C^2$, so $A_5$ has an action on $\Proj^1$ with abelian stabilizers. There is a similar statement for PNC central extensions.

Secondly, due to the Chevalley-Shepard-Todd theorem, if the point stabilizers contain pseudoreflections, it is not actually necessary for them to be abelian in order for the resulting quotient singularities to be (at worst) abelian. For example, although neither $A_5$ nor $PSL(2,7)$ are PNC, if $V$ is a three-dimensional faithful irreducible representation of $PSL(2,7)$ over $\C$, then $\Proj(V)/PSL(2,7)$ has only abelian singularities, while if $W$ is a three-dimensional faithful irreducible representation of $A_5$, then $\Proj(W)/A_5$ is actually {\em smooth}. Details are found in \cite[Section 1.2]{blumsmith}. This prompts us to ask:

\begin{question}
For which nonabelian finite simple groups $G$ is there a faithful representation $V$ such that $\Proj(V)/G$ has at worst abelian singularities?
\end{question}

Of course, one can broaden this question beyond simple groups.

\section{Acknowledgement}

The authors wish to thank Frieder Ladisch for the statement and proof of lemma \ref{ladischlemma}.

The second author has been funded by
the Russian Academic Excellence Project `5-100' and aknowledges support by a Simons travel grant and by the EPSRC program grant EP/M024830.

With the exception of sections \ref{sec:permutations} and \ref{sec:questions}, this paper is a revision of work that originally appeared as part of the first author's doctoral thesis under the supervision of the second author and Yuri Tschinkel at the Courant Institute of Mathematical Sciences at NYU.

\end{document}